\documentclass[twoside]{article} 
\usepackage{graphicx}
\date{} 
\title{The asymptotic expansion of the Bernoulli polynomials of the second kind}
\author{\sc R. B.\ Paris \\
{\em Division of Computing and Mathematics,} \\
{\em Abertay University, Dundee DD1 1HG, UK}}
\begin{document}
\def\f#1#2{\mbox{${\textstyle \frac{#1}{#2}}$}}
\def\dfrac#1#2{\displaystyle{\frac{#1}{#2}}}
\def\boldal{\mbox{\boldmath $\alpha$}}
\newcommand{\bee}{\begin{equation}}
\newcommand{\ee}{\end{equation}}
\newcommand{\sa}{\sigma}
\newcommand{\ka}{\kappa}
\newcommand{\al}{\alpha}
\newcommand{\la}{\lambda}
\newcommand{\ga}{\gamma}
\newcommand{\eps}{\epsilon}
\newcommand{\om}{\omega}
\newcommand{\fr}{\frac{1}{2}}
\newcommand{\fs}{\f{1}{2}}
\newcommand{\g}{\Gamma}
\newcommand{\br}{\biggr}
\newcommand{\bl}{\biggl}
\newcommand{\ra}{\rightarrow}
\newcommand{\gtwid}{\raisebox{-.8ex}{\mbox{$\stackrel{\textstyle >}{\sim}$}}}
\newcommand{\ltwid}{\raisebox{-.8ex}{\mbox{$\stackrel{\textstyle <}{\sim}$}}}
\renewcommand{\topfraction}{0.9}
\renewcommand{\bottomfraction}{0.9}
\renewcommand{\textfraction}{0.05}
\newcommand{\mcol}{\multicolumn}
\date{}
\maketitle
\pagestyle{myheadings}
\markboth{\hfill \sc R. B.\ Paris  \hfill}
{\hfill \sc Asymptotics of the Bernoulli polynomial\hfill}
\begin{abstract}
We consider the Bernoulli polynomials of the second kind, which can be related to the generalised Bernoulli polynomials $B_n^{(n)}(z)$. The asymptotic expansions of the scaled polynomials $B_n^{(n)}(nz)$ are obtained as $n\to\infty$ when (i) $z$ is real and (ii) $z$ is complex bounded away from $[0,1]$. These results complement recent work of \v{S}tampach [{\it J. Approx. Theory}, {\bf 262} (2021) 105517].
Numerical results are presented to illustrate the accuracy of the different expansions obtained.
\vspace{0.3cm}

\vspace{0.1cm}
 
\end{abstract}

\vspace{0.3cm}

\noindent $\,$\hrulefill $\,$

\vspace{0.3cm}
\begin{center}
{\bf 1.\ Introduction}
\end{center}
\setcounter{section}{1}
\setcounter{equation}{0}
\renewcommand{\theequation}{\arabic{section}.\arabic{equation}}
The generalised Bernoulli polynomials of order $\mu$ are defined by the generating function
\bee\label{e12}
\bl(\frac{t}{e^t-1}\br)^{\!\mu} e^{zt}=\sum_{n=0}^\infty B_n^{(\mu)}(z)\,\frac{t^n}{n!} \qquad (|t|<2\pi)
\ee
and satisfy the well-known property
\bee\label{e13}
B_n^{(\mu)}(-z)=(-)^n B_n^{(\mu)}(\mu+z).
\ee
The Bernoulli polynomials of the second kind $b_n(z)$ are defined by
\[\frac{t(1+t)^z}{\log (1+t)}=\sum_{n=0}^\infty b_n(z) \frac{t^n}{n!} \qquad (|t|<1)\]
and are related to the generalised Bernoulli polynomials by
\bee\label{e11}
b_n(z)=B_n^{(n)}(z+1).
\ee
Consequently, it will be sufficient in what follows to consider the polynomials $B_n^{(n)}(z)$ instead of $b_n(z)$ on account of the relation (\ref{e11}). 
The first few polynomials $B_n^{(n)}(z)$ are
\begin{eqnarray*}
B_0^{(0)}(z)&=&1,\quad B_1^{(1)}(z)=\frac{1}{2} (-1+2z),\quad B_2^{(2)}(z)=\frac{1}{6}(5-12z+6z^2),\\
B_3^{(3)}(z)&=&  \frac{1}{4}(-9+24z-18z^2+4z^3),\\
B_4^{(4)}(z)&=&\frac{1}{30}(251-720z+66z^2-240z^3+30z^4),\\
B_5^{(5)}(z)&=& \frac{1}{12}(-475+1449z-1500z^2+700z^3-150z^4+12z^5), \ldots \ .
\end{eqnarray*}

The asymptotic expansion of $B_n^{(\mu)}(z)$ for large $\mu$ (fixed $n$) and large $n$ (fixed $\mu$) has been investigated in \cite{LT1,LT2}.
The asymptotic behaviour when $\mu=n$ for large $n$ and the structure of the zeros have been studied recently by \v{S}tampach \cite{FS}. It was shown that, unlike the zeros of the Bernoulli polynomials of the first kind, the zeros of the Bernoulli polynomials of the second kind are all real, simple and confined to the interval $z\in[0,n]$. In addition, he showed that the zeros are distributed symmetrically about the point $z=n/2$ and interlace with the integers $1, 2, \ldots , n$. 

\v{S}tampach  determined the behaviour of the zeros near the endpoints 0 and $n$ and also of those zeros located at a fixed distance from the point $\al n$, where $\al\in (0,1)$, by obtaining the expansion of $B_n^{(n)}(z+\al n)$ for $z\in {\bf C}$. Particular attention was paid to the case $\al=\fs$, namely for those zeros situated in the neighbourhood of the mid-point. When considering the scaled polynomials $B_n^{(n)}(nz)$, the oscillatory region is $z\in[0,1]$; \v{S}tampach also determined the leading asymptotic form in the non-oscillatory (zero-free) region $z\in {\bf C}\backslash[0,1]$.

In this article, we complement the analysis of \v{S}tampach \cite{FS} by determining the asymptotic expansion of $B_n^{(n)}(nz)$
for $n\to\infty$ and $z\in{\bf C}$. From (\ref{e13}) and the conjugacy property, we have
\bee\label{e14}
B_n^{(n)}(nz)=(-)^n B_n^{(n)}(n(1-z)),\qquad B_n^{(n)}(n{\overline z})={\overline {B_n^{(n)}(nz)}},\
\ee
so that it is sufficient to consider the expansion in the quadrant $\Re (z)\geq\fs$, $\Im (z)\geq 0$.
Unlike \v{S}tampach, who used different representations for the different cases, we use a single integral representation combined with the method of steepest descents.
\vspace{0.6cm}

\begin{center}
{\bf 2.\ The asymptotic expansion of $B_n^{(n)}(nx)$ for $x>0$}
\end{center}
\setcounter{section}{2}
\setcounter{equation}{0}
\renewcommand{\theequation}{\arabic{section}.\arabic{equation}}
We consider the scaled generalised Bernoulli polynomials $B_n^{(n)}(nz)$ for $n\to\infty$, for which the oscillatory region is $z\in[0,1]$.  From (\ref{e12}), we obtain the integral representation 
\bee\label{e21}
B_n^{(n)}(nz)=\frac{n!}{2\pi i} \oint \frac{e^{nzs}}{(e^s-1)^n\, s}\,ds=\frac{n!}{2\pi i} \oint e^{-n\psi(s)}\,\frac{ds}{s},
\ee
where 
\[\psi(s)\equiv \psi(s,z):=\log (e^s-1)-sz\]
and the integration path is a closed circuit described in the positive sense surrounding the origin but excluding the poles $\pm2\pi i$. Saddle points of $\psi(s)$ occur at $\psi'(s)=0$; that is when $e^s/(e^s-1)-z=0$. Taking account of the multi-valued nature of the logarithm, we have the saddles $s_k$ given by
\bee\label{e22}
s_k=\log \frac{z}{z-1}+2\pi ik \qquad (k=0, \pm 1, \pm 2, \ldots ).
\ee
In this section we deal with the case when $z$ $(=x)$ is real; from (\ref{e14}), it is sufficient to consider only the case $x>0$ (or more precisely $x\geq\fs$). 
\medskip

\noindent{\bf 2.1\ \ The case $x>1$}
\medskip

\noindent
When $x>1$ it is seen that $s_0=\log (x/(x-1))$ is situated on the positive real axis and that $\Im \psi(s_0)=0$. The steepest descent path through $s_0$ is given by $\Im \psi(s)=0$; that is, with $s=\xi+i\eta$,
\bee\label{e22a}
\eta(1-x)+\phi=0,\qquad \phi:=\arctan \bl(\frac{\sin \eta}{e^\xi-\cos \eta}\br).
\ee
When $\xi\to-\infty$, we have $\phi\sim \pm\pi-\eta$; consequently, the steepest descent path through $s_0$ passes to $-\infty$ with $\eta\sim\pm\pi/x$; see Fig.~1(a).

\begin{figure}[th]
	\begin{center}	{\tiny($a$)}\ \includegraphics[width=0.4\textwidth]{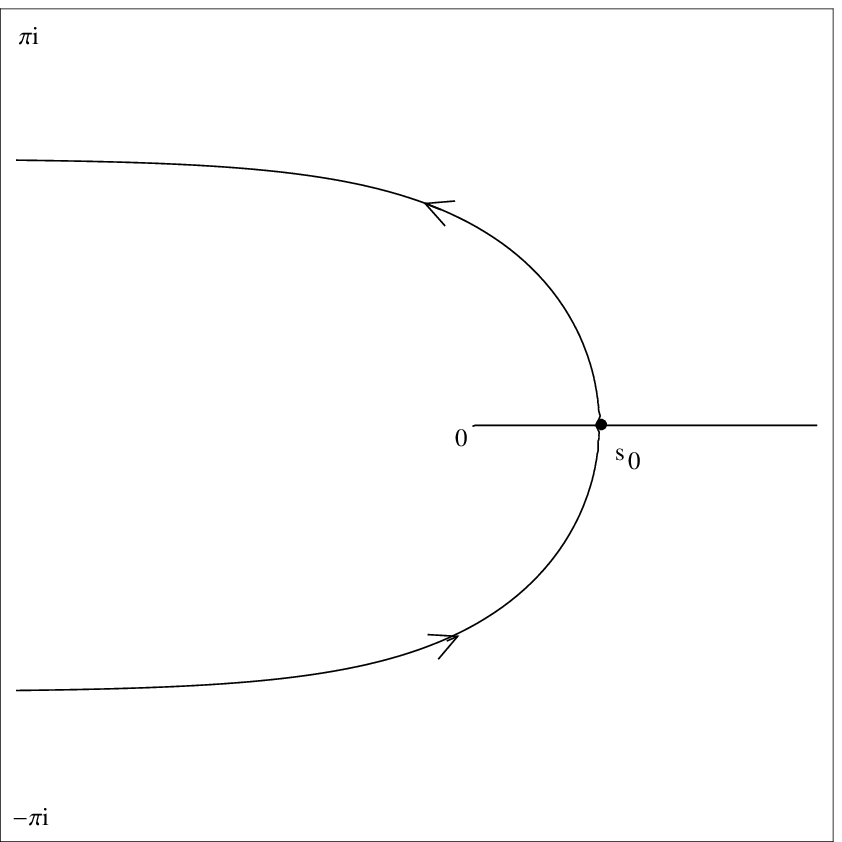}
	\qquad
	{\tiny($b$)}\ \includegraphics[width=0.368\textwidth]{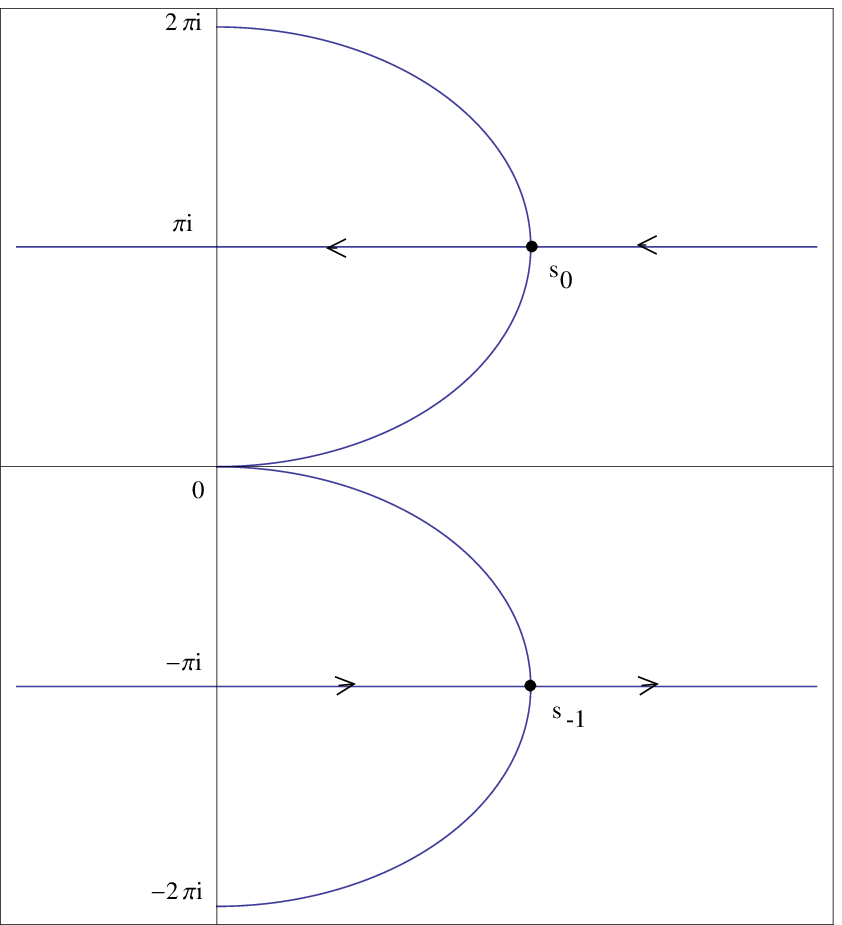} 
	
\caption{\small{Paths of steepest descent when $x>0$ (the arrows indicate the direction of integration): ($a$) $x>1$ when the contributory saddle is $s_0$ and the steepest ascent path is the positive real axis and ($b$) $0<x<1$ when the contributory saddles are $s_0$ and $s_{-1}$ with the steepest ascent paths passing from the origin to $\pm2\pi i$.}}
\end{center}
\end{figure}

Then, from (\ref{e21}), we have
\bee\label{e23}
B_n^{(n)}(nx)
=\frac{n!}{2\pi i}\int_C e^{-n\psi(s)} \frac{ds}{s}
=\frac{n!}{2\pi i}\, e^{-n\psi(s_0)} \int_{-\infty}^\infty e^{-\frac{1}{2}nw^2} \frac{1}{s}\frac{ds}{dw}\,dw,
\ee
where we have introduced the new variable $\psi(s)-\psi(s_0)=\fs w^2$ and $C$ is the steepest descent path.
If we let
\[s_0=\log\,h,\qquad h=\frac{x}{x-1},\]
we obtain the expansion
\[\psi(s)-\psi(s_0)=\frac{\psi''(s_0)}{2!} (s-s_0)^2+\frac{\psi'''(s_0)}{3!} (s-s_0)^3+\frac{\psi^{iv}(s_0)}{4!} (s-s_0)^4+\cdots\]
\[=-\frac{h}{2! (h-1)^2} (s-s_0)^2+\frac{h(h+1)}{3! (h-1)^3} (s-s_0)^3-\frac{h(1+4h+h^2)}{4! (h-1)^4} (s-s_0)^4+\ldots =\fs w^2,\]
which upon inversion yields
\[s-s_0=i\frac{h-1}{h^{1/2}}-\frac{(h^2-1)}{6h} w^2+i\frac{(h-1)(-1+h-h^2)}{36 h^{3/2}} w^3+\ldots\, .\]
Then
\bee\label{e24}
\frac{1}{s} \frac{ds}{dw} \stackrel{e}{=} \frac{i(h-1)}{s_0 h^{1/2}}\bl\{1+A_1(h,s_0)w^2+A_2(h,s_0)w^4+A_3(h,s_0) w^6+\cdots \br\},
\ee
where $\stackrel{e}{=}$ signifies the inclusion of {\em only the even powers of $w$}, since odd powers will not enter into the calculations. With the help of {\it Mathematica} the first three coefficients\footnote{The first two coefficients $A_j(h,s_0)$ can be obtained alternatively from the standard representation of the saddle-point coefficients; see, for example, \cite[pp.~13--14]{PHad}.} $A_j(h,\la)$ are given by
\[A_1(h,\la)=\frac{1}{12h}\bl\{-(1-h+h^2)+\frac{6(h^2-1)}{\la}-\frac{12(h-1)^2}{\la^2}\br\},\]
\[A_2(h,\la)=\frac{1}{864h^2}\bl\{(1-h+h^2)^2-\frac{12}{\la}(h^2-1)(3-5h+3h^2)+\frac{120}{\la^2}(h-1)^2(2+h+2h^2)\]
\[\hspace{4cm}-\frac{720}{\la^3} (h+1)(h-1)^3+\frac{864}{\la^4} (h-1)^4\br\},\]
\[A_3(h,\la)=\frac{1}{777600h^3}\bl\{\Upsilon(h)+\frac{90}{\la} (h^2-1)(1-h+h^2)(5-9h+5h^2)-\frac{1260}{\la^2}(h-1)^2(13-8h-3h^2-8h^3+13h^4)\]
\[\hspace{2cm}+\frac{15120}{\la^3}(h-1)^3(h+1)(8-5h+8h^2)-\frac{453600}{\la^4}(h-1)^4(1+h+h^2)\]
\bee\label{e2coeff}
\hspace{3cm}+\frac{907200}{\la^5}(h-1)^5 (h+1)-\frac{777600}{\la^6} (h-1)^6\br\},
\ee
where
\[\Upsilon(h)=139-417h+402h^2-109h^3+402h^4-417h^5+139h^6.\]
Higher order coefficients can be calculated in this manner when $x$ is given a specific value; see Section 4 for an example.

Upon noting that
\[e^{-n\psi(s_0)}=(x-1)^n \bl(\frac{x}{x-1}\br)^{\!nx},\]
we find, on substituting the expansion (\ref{e24}) into (\ref{e23}),
\[B_n^{(n)}(nx)\sim \frac{n!}{2\pi}\,\frac{(x-1)^{n-1}}{s_0} \bl(\frac{x}{x-1}\br)^{\!nx-1/2} \int_{-\infty}^\infty
e^{-\frac{1}{2}nw^2} \sum_{k=0}^\infty A_k(h,s_0) w^k\,dw.\] 
Evaluation of the integrals in terms of gamma functions then leads to the following result:
\newtheorem{theorem}{Theorem}
\begin{theorem}$\!\!\!.$\ \ As $n\to\infty$ when $x>1$, we have the expansion
\bee\label{e25}
B_n^{(n)}(nx)\sim \frac{n!}{\sqrt{2\pi n}}\,\frac{(x-1)^{n-1}}{\log\,(x/(x-1))} \bl(\frac{x}{x-1}\br)^{\!nx-1/2} \sum_{k=0}^\infty \frac{2^k(\fs)_k }{n^k}\,A_k(h,s_0),
\ee
where $h=x/(x-1)$ and the coefficients $A_k(h,\la)$ are given in (\ref{e2coeff}) with $\la=s_0=\log\,h$.
\end{theorem}

The leading term of (\ref{e25}) can be shown to agree with \cite[Eq.~(40)]{FS} (when the variable $z$ therein vanishes).
\medskip

\noindent{\bf 2.2\ \ The case $0<x<1$}
\medskip

\noindent
When $0<x<1$, it is sufficient to restrict attention to the interval $\fs\leq x<1$ by (\ref{e14}).
The contributory saddles are $s_0=\log\,h+\pi i$ and $s_{-1}=\log\,h-\pi i$, where now $h=x/(1-x)$.
We have $\Im (\psi(s_0))=\pi(1-x)$ and, with $s=\xi+i\eta$, $\Im (\psi(s))=\eta(1-x)+\phi$, where $\phi$ is defined in (\ref{e22a}). Then it is seen that when $\eta=\pi$, $\phi=0$ and $\Im (\psi(s))=\Im (\psi(s_0))$. The steepest descent path through $s_0$ is thus the line $\eta=\pi$; a similar argument shows that $\eta=-\pi$ is the steepest descent path through $s_{-1}$. Consequently the contour $C$ can be taken as the horizontal paths $(-\infty-\pi i, \infty-\pi i)$ and $(\infty+\pi i,-\infty+\pi i)$; see Fig.~1(b).

Consider the lower saddle $s_{-1}$. As in Section 2.1, we have
\[\frac{1}{s} \frac{ds}{dw} \stackrel{e}{=}  \frac{h+1}{s_{-1} h^{1/2}} \bl\{1+A_1(-h,s_{-1}) w^2+A_2(-h,s_{-1})w^4+A_3(-h,s_{-1}) w^6+\cdots \br\},\]
where the coefficients $A_k(-h,s_{-1})$ are given in (\ref{e2coeff}) with $\la=s_{-1}$ and
\bee\label{e25a}
s_{-1}=\log\,h-\pi i=Le^{-i\om},\qquad L:=\{\log^2h+\pi^2\}^{1/2},\ \ \om:=\arctan\,(\pi/\log\,h).
\ee
Then, since
\[e^{-n\psi(s_{-1})}=(-)^n (1-x)^n \bl(\frac{x}{1-x}\br)^{\!nx} e^{-\pi inx},\]
we find the contribution
\[\frac{e^{-n\psi(s_{-1})}}{2\pi i} \int_{-\infty}^\infty e^{-\frac{1}{2}nw^2}\frac{1}{s} \frac{ds}{dw}\,dw\hspace{7cm}\]
\[\sim \frac{(-)^n}{2\pi L} (1-x)^{n-1} \bl(\frac{x}{1-x}\br)^{\!nx-1/2} e^{-i\Theta(x)} \int_{-\infty}^\infty e^{-\frac{1}{2}nw^2} \sum_{k=0}^\infty A_k(-h,s_{-1}) w^k\,dw\] 
\[=\frac{(-)^n}{\sqrt{2\pi n}\,L}(1-x)^{n-1} \bl(\frac{x}{1-x}\br)^{\!nx-1/2} e^{-i\Theta(x)} \sum_{k=0}^\infty \frac{2^k(\fs)_k}{n^k}\,A_k(-h,s_{-1}),\hspace{1.1cm}\]
where
\[\Theta(x):=\pi nx-\om+\fs\pi.\]

The contribution from the upper saddle $s_0$ is given by the conjugate expression. We observe that due to the presence of $\la=s_{-1}$ in the coefficients (\ref{e2coeff}) in the case $0<x<1$, the coefficients $A_k(-h,s_{-1})$ separate into real and imaginary parts. For example, in the case of $A_1(-h,s_{-1})$ we find
\[\Re (A_i(-h,s_{-1}))=\frac{1}{12h}\bl\{1+h+h^2-\frac{6(h^2-1)}{L} \cos \om+\frac{12(h+1)^2}{L^2} \cos 2\om\br\},\]
\[\Im (A_1(-h,s_{-1}))=\frac{1}{12h}\bl\{-\frac{6(h^2-1)}{L} \sin \om+\frac{12(h+1)^2}{L^2} \sin 2\om\br\}.\hspace{1.9cm}\]
Hence, upon some routine algebra, we obtain the following result:
\begin{theorem}$\!\!\!.$\ \ Let $h=x/(1-x)$ and $\Theta(x)=\pi nx-\om+\fs\pi$. Then as $n\to\infty$ when $0<x<1$, we have the expansion
\[B_n^{(n)}(nx)\sim\frac{2(-)^n n!}{\sqrt{2\pi n}}\,\frac{(1-x)^{n-1}}{\{\log^2h+\pi^2\}^{1/2}} \bl(\frac{x}{1-x}\br)^{\!nx-1/2}\hspace{5cm}\]
\bee\label{e26}
\hspace{1cm}\times\bl\{\cos \Theta(x) \sum_{k=0}^\infty \frac{2^k(\fs)_k}{n^k}\,\Re (A_k(-h,s_{-1}))+\sin \Theta(x) \sum_{k=1}^\infty \frac{2^k(\fs)_k}{n^k}\,\Im (A_k(-h,s_{-1}))\br\}.
\ee
The coefficients $A_k(-h,\la)$ are given by (\ref{e2coeff}) with the quantity $\la=s_{-1}=Le^{-i\om}$, where $L$ and $\om$ are defined in (\ref{e25a}). 
\end{theorem}
\medskip

\noindent{\bf 2.3\ \ The case $x=\fs$}
\medskip

\noindent
When $x=\fs$, we have $h=1$, $s_{-1}=-\pi i$, $L=\pi$, $\om=\fs\pi$, $\Theta(x)=\fs\pi n$ and $\Im (A_k(-1,\la))=0$. Then we obtain from (\ref{e26}) the expansion
\bee\label{e201}
B_n^{(n)}(\fs n)\sim\frac{2^{2-n} n!}{\sqrt{2\pi n}}\,\frac{\cos \fs\pi n}{\pi} \sum_{k=0}^\infty \frac{(-2)^k(\fs)_k}{n^k}\,C_k\qquad(n\to\infty),
\ee
where $C_k=(-)^k\Re (A_k(-1,-\pi i))$ and $n$ is even; when $n$ is odd $B_n^{(n)}(\fs n)\equiv 0$. In this case it is possible to carry the algebra described in Section 2.1 further to find the coefficients
\[C_0=1,\quad C_1=\frac{1}{4\pi^2}(16 - \pi^2), \quad C_2=\frac{1}{96\pi^4}(1536 - 160 \pi^2 + \pi^4),\]
\[C_3=\frac{1}{5760\pi^6}(368640 - 53760 \pi^2 + 1456 \pi^4 + 15 \pi^6),\]
\[C_4=\frac{1}{645120 \pi^8}(165150720 - 30965760 \pi^2 + 1483776 \pi^4 - 3904 \pi^6 - 63 \pi^8),\]
\[C_5=\frac{1}{38707200 \pi^{10}}(39636172800 - 9083289600 \pi^2 + 624476160 \pi^4 -
 10081280 \pi^6 - 92048 \pi^8 - 1995 \pi^{10}). \]
The first four coefficients $C_k$ ($k\leq 3$) follow from (\ref{e2coeff}) with $\la=-\pi i$. The expansion (\ref{e201}) agrees with that given in \cite[Eq.~(50)]{FS} when $z=0$.

\vspace{0.6cm}
\newpage
\begin{center}
{\bf 3.\ The asymptotic expansion of $B_n^{(n)}(nz)$ for complex $z$}
\end{center}
\setcounter{section}{3}
\setcounter{equation}{0}
\renewcommand{\theequation}{\arabic{section}.\arabic{equation}}
As stated in Section 1, when $z=x+iy$ is complex it is sufficient to restrict attention to the quadrant $x\geq \fs$ and $y>0$, so that $0<\arg\,z<\fs\pi$. The contributory saddles in this case are found to be $s_0$ and $s_1$, where we note that
\bee\label{e31}
\psi(s_1)=\psi(s_0+2\pi i)=\psi(s_0)+2\pi i(1-z)
\ee
and some straightforward algebra shows that $-\fs\pi\leq \arg\,s_0<0$ when $x\geq\fs$, $y>0$.
Fig.~2 shows typical paths of steepest descent when ($a$) $\fs\leq x<1$, ($b$) $x=1$ and ($c$) $x>1$.

\begin{figure}[th]
	\begin{center}	{\tiny($a$)}\ \includegraphics[width=0.26\textwidth]{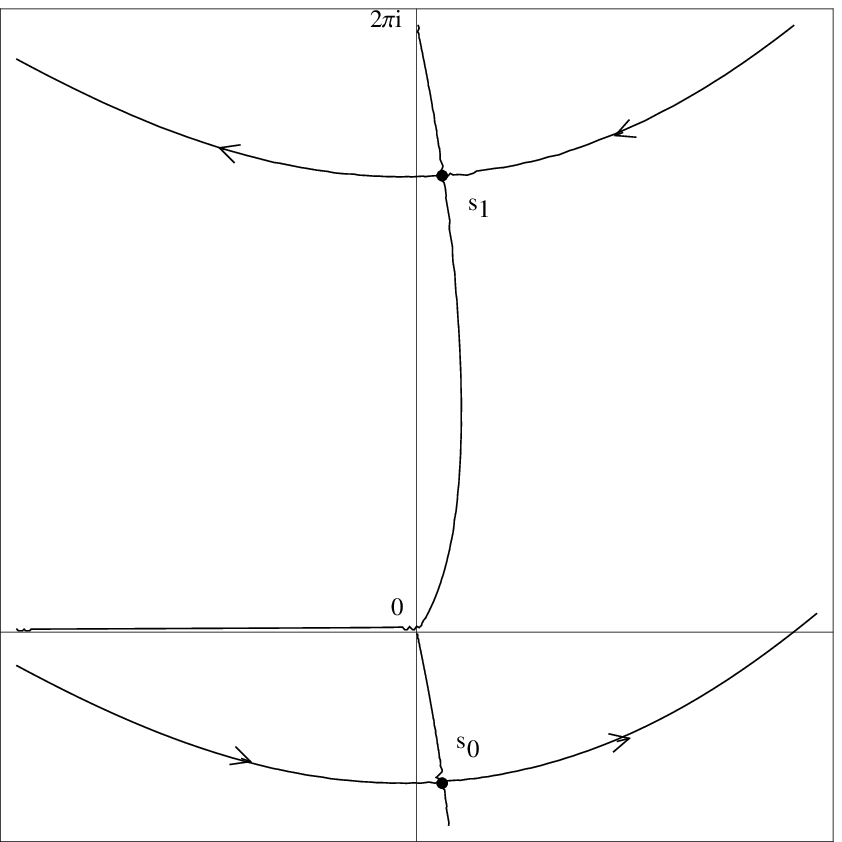}
	\qquad
	{\tiny($b$)}\ \includegraphics[width=0.26\textwidth]{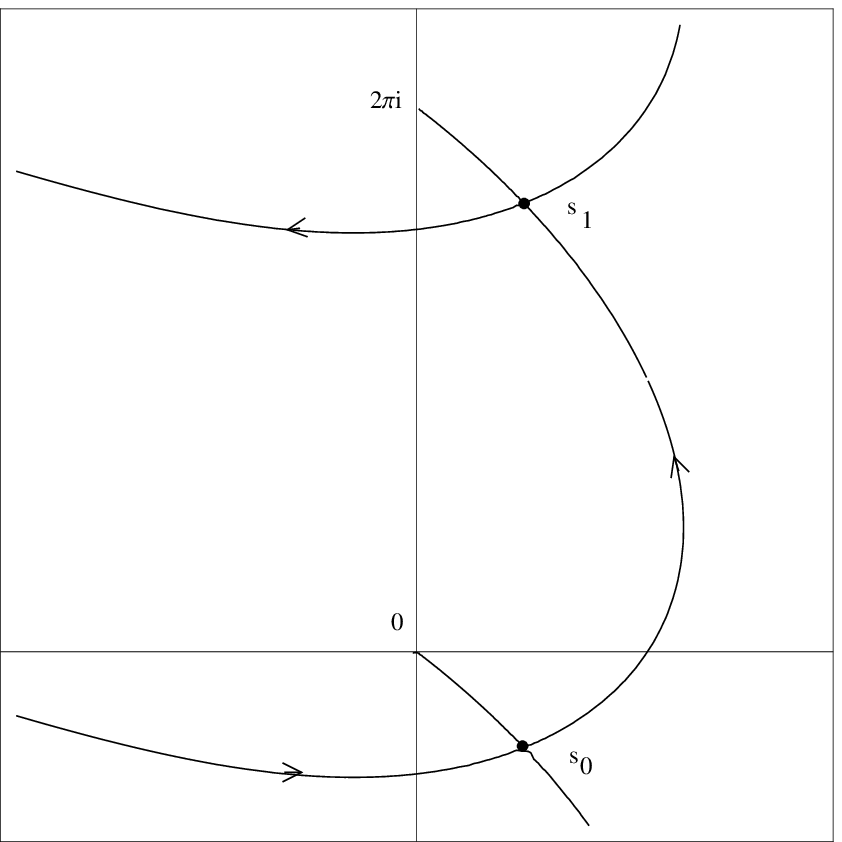}
	\qquad 
	{\tiny($c$)}\ \includegraphics[width=0.26\textwidth]{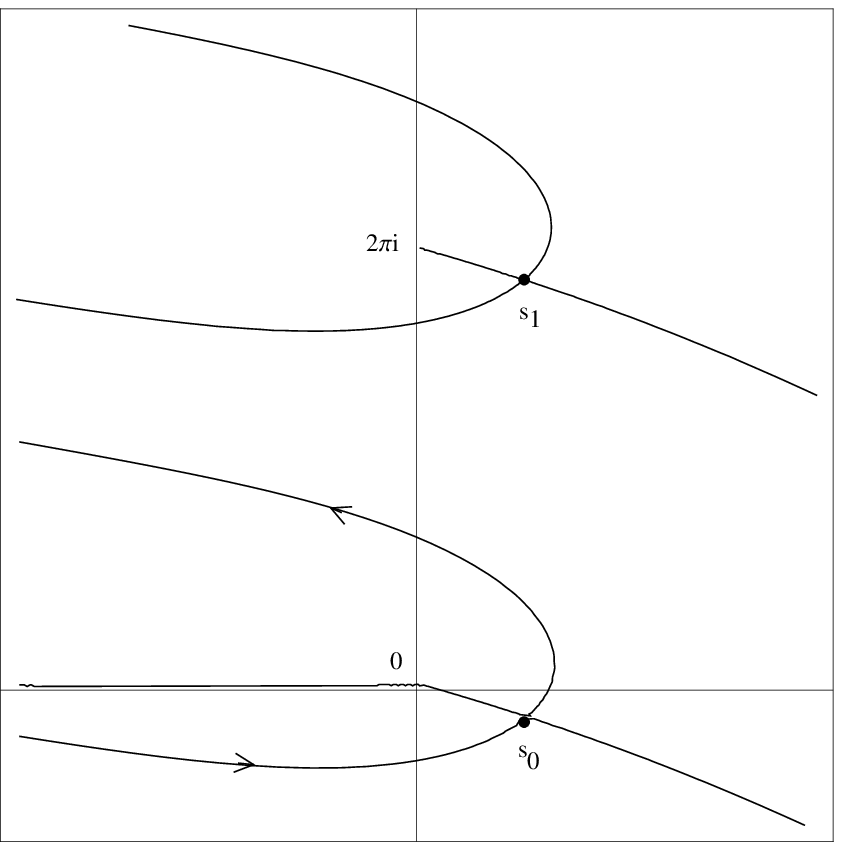}
\caption{\small{Paths of steepest descent when $z=x+iy$, $y>0$: ($a$) $\fs\leq x<1$, ($b$) $x=1$ and ($c$)  $x>1$. The arrows indicate the direction of integration.}}
\end{center}
\end{figure}

When $\fs\leq x<1$, $y>0$ the steepest descent paths that pass to infinity in $\Re (s)<0$ have slope $\simeq\pi-\arctan (y/x)$, whereas those in $\Re (s)>0$ have slope $\simeq \arctan (y/(1-x))$. The closed path in (\ref{e21}) can be reconciled with the steepest paths through $s_0$ and $s_1$ in the directions indicated in Fig.~2(a). The direction of integration through $s_0$ is given by $\fs\pi-\fs\arg\,\psi''(s_0)$, where $\psi''(s_0)=z(1-z)$.  When $x=1$, the steepest path through $s_0$ connects\footnote{From(\ref{e31}), we have $\Im (\psi(s_0)-\psi(s_1))=2\pi (x-1)$ which vanishes when $x=1$.} with the saddle $s_1$ (a Stokes phenomenon); again the closed contour can be reconciled with the steepest path thorough $s_0$ and half of the path emanating from $s_1$; see Fig.~2(b). Finally, when $x>1$ the closed contour in (\ref{e21}) is reconcilable with the steepest path through $s_0$ alone; see Fig.~2(c). 

Then, with $h=z/(z-1)$ and $s_0=\log\,h$, the contribution from the saddle $s_0$ is given by the formal asymptotic sum
\bee\label{e32}
S_0(z):=\frac{n!}{\sqrt{2\pi n}}\,\frac{(z-1)^{n-1}}{\log\,(z/(z-1))} \bl(\frac{z}{z-1}\br)^{\!nz-1/2} \sum_{k=0}^\infty \frac{2^k (\fs)_k}{n^k}\,A_k(h,s_0).
\ee
From (\ref{e31}), we have
\[e^{-n\psi(s_1)}=e^{-n\psi(s_0)}e^{-2\pi ni(1-z)}=e^{-n\psi(s_0)} O(e^{-2\pi ny}).\]
Thus, for fixed $y>0$, the contribution from the saddle $s_1$ is exponentially small as $n\to\infty$. This contribution is given by the formal asymptotic sum
\bee\label{e33}
S_1(z):=-\frac{n!}{\sqrt{2\pi n}}\,\frac{e^{-2\pi ni(1-z)}(z-1)^{n-1}}{(\log\,(z/(z-1))+2\pi i)} \bl(\frac{z}{z-1}\br)^{\!nz-1/2} \sum_{k=0}^\infty \frac{2^k (\fs)_k}{n^k}\,A_k(h,s_1)
\ee
when $\fs\leq x<1$, $y>0$. We do not consider here the exponentially small contribution when $x=1$.  

Collecting together (\ref{e32}) and (\ref{e33}), we have the result:
\begin{theorem}$\!\!\!.$\ \ Let $h=z/(z-1)$, $s_0=\log\,h$ and $s_1=s_0+2\pi i$. Then as $n\to\infty$ and $z$ bounded away from $[0,1]$ we have the expansions
\bee\label{e34}
B_n^{(n)}(nz)\sim \left\{\begin{array}{ll} S_0(z)+S_1(z) & (\fs\leq x<1,\ y>0)\\
\\
S_0(z) & (x\geq 1,\ y>0), \end{array}\right.
\ee
where $S_0(z)$ and $S_1(z)$ are defined by the formal asymptotic sums (\ref{e32}) and (\ref{e33}). The expansion $S_1(z)$ is exponentially smaller than $S_0(z)$  by the factor $O(e^{-2\pi ny})$ and may, in most applications, be neglected. The first few coefficients $A_k(h,\la)$ associated with $S_0(z)$ and $S_1(z)$ with $\la=s_0$ and $\la=s_1$ respectively, are given in (\ref{e2coeff}).
\end{theorem}
The leading form of the expansion $S_0(z)$ was given by \v{S}tampach \cite[Eq.~(61)]{FS}.
\medskip

\noindent{\bf Remark 1.}\ \ We observe that the expansion $S_0(z)$ when $x>1$, $y>0$ has the same form as that given in (\ref{e25}) for $x>1$, $y=0$. Hence, it follows that the domain of validity of the second expansion in (\ref{e34}) 
can also include the domain $x>1$, $y\geq 0$. 

\vspace{0.6cm}

\begin{center}
{\bf 4.\ Numerical results and concluding remarks}
\end{center}
\setcounter{section}{4}
\setcounter{equation}{0}
\renewcommand{\theequation}{\arabic{section}.\arabic{equation}}
To illustrate the accuracy of the expansions in Theorems 1--3, we show the values of the absolute relative error in the computation of $B_n^{(n)}(nz)$ using the asymptotic expansions in (\ref{e25}), (\ref{e26}) and (\ref{e34}) for different $n$ and $z$ as a function of the truncation index $k$. The value of $B_n^{(n)}(nz)$ was computed in {\it Mathematica} by the command NorlundB[n, n, nz]. In Table 1 we show\footnote{In the tables we have adopted the convention of writing $x(y)$ to represent $x\times 10^{y}$.} the absolute relative errors for complex $z$ using the expansion $S_0(z)$ in Theorem 3. Table 2 shows the errors
when $z$ is real on the interval $[\fs,1)$ using the expansion in Theorem 2.
\begin{table}[h]
\caption{\footnotesize{The absolute relative error in the computation of $B_n^{(n)}(nz)$ for different $z$ when $n=20$ as a function of the truncation index $k$.}}
\begin{center}
\begin{tabular}{|c|c|c|c|c|c|}
\hline
&&&&&\\[-0.3cm]
\mcol{1}{|c|}{$k$} & \mcol{1}{c|}{$z=2$} & \mcol{1}{c|}{$z=2+i$} & \mcol{1}{c|}{$z=1+i$} & \mcol{1}{c|}{$z=0.75+i$}
& \mcol{1}{c|}{$z=0.50+i$}\\
\hline
&&&&&\\[-0.3cm]
0 & $4.193(-3)$ & $4.169(-3)$ & $4.173(-3)$ & $4.208(-3)$ & $4.227(-3)$\\
1 & $7.780(-6)$ & $9.459(-6)$ & $9.250(-6)$ & $7.304(-6)$ & $6.123(-6)$\\
2 & $3.449(-7)$ & $3.323(-7)$ & $3.413(-7)$ & $3.378(-7)$ & $3.344(-7)$\\
3 & $1.051(-9)$ & $1.779(-9)$ & $1.630(-9)$ & $2.588(-9)$ & $2.843(-9)$\\  
\hline
\end{tabular}
\end{center}
\end{table}
\begin{table}[h]
\caption{\footnotesize{The absolute relative error in the computation of $B_n^{(n)}(nx)$ for different values of $x$ and $n$ as a function of the truncation index $k$.}}
\begin{center}
\begin{tabular}{|c|cc|cc|}
\hline
&&&&\\[-0.3cm]
\mcol{1}{|c|}{} & \mcol{2}{c|}{$n=20$} & \mcol{2}{c|}{$n=40$}\\
\mcol{1}{|c|}{$k$} & \mcol{1}{c}{$x=0.50$} & \mcol{1}{c|}{$x=0.75$} & \mcol{1}{c}{$x=0.50$} & \mcol{1}{c|}{$x=0.75$}\\
\hline
&&&&\\[-0.3cm]
0 & $7.719(-3)$ & $8.725(-3)$ & $3.871(-3)$ & $4.370(-3)$ \\
1 & $4.578(-5)$ & $2.597(-5)$ & $1.116(-5)$ & $5.556(-6)$ \\
2 & $1.909(-6)$ & $7.264(-6)$ & $2.402(-7)$ & $9.012(-7)$ \\
3 & $3.620(-8)$ & $8.782(-8)$ & $3.094(-9)$ & $8.062(-8)$ \\
\hline
\end{tabular}
\end{center}
\end{table}

In order to verify the assertion made in Theorem 3 concerning the appearance of the exponentially small expansion $S_1(z)$ when $\fs\leq x<1$, $y>0$, it is necessary to select values of $n$ and $\Im (z)$ not too large. To detect the presence of $S_1(z)$ we have to optimally truncate the dominant expansion $S_0(z)$ at, or near, its least term in modulus. We choose to work with the value $n=10$, for which the three explicit representations of the coefficients $A_k(h,s_0)$ in (\ref{e2coeff})
are insufficient to achieve optimal truncation. The procedure described in Section 2.1 of expansion and inversion can be applied in specific cases where the value of $z$ is specified to produce numerical values of the coefficients for high $k$-values. These coefficients are shown in Table 3 for the particular case $z=2/3+i/4$, for which it is found that optimal truncation of $S_0(z)$ when $n=10$ occurs at $k=10$.
\begin{table}[h]
\caption{\footnotesize{The coefficients $A_k(h,s_0)$ for $1\leq k\leq 10$  for $z=2/3+i/4$ (with $A_0(h,s_0)=1$).}}
\begin{center}
\begin{tabular}{|c|l|}
\hline
&\\[-0.3cm]
\mcol{1}{|c|}{$k$} & \mcol{1}{c|}{$A_k(h,s_0)$}\\
\hline
&\\[-0.3cm]
1 & $-1.0029378942(-01) -1.8804724469i(-02)$\\ 
2 & $-3.7372334426(-03) -5.5650719166i(-04)$\\
3 & $+1.8095948417(-05) +1.5684946154i(-04)$\\ 
4 & $+5.9175620462(-05) +1.3608152444i(-04)$\\ 
5 & $+5.6624929259(-06) +5.2629558202i(-06)$\\ 
6 & $+3.2408350155(-03) -2.4032813980i(-06)$\\
7 & $+8.6041310199(-08) -2.5286915962i(-07)$\\
8 & $-1.0224648657(-07) -8.6048696324i(-08)$\\
9 & $-8.4341941837(-09) -3.2178913880i(-10)$\\
10& $+5.6624929259(-06) +5.2629558202i(-06)$\\
\hline
\end{tabular}
\end{center}
\end{table}

In Table 4, we show the values of $B_n^{(n)}(nz)-S_0(z)$, where the dominant expansion $S_0(z)$ is optimally truncated,
compared with the exponentially small expansion $S_1(z)$ (with $k\leq 3$) when $n=10$ and $z=x+i/4$. It is seen that when $x<1$ there is good agreement between these two values, thereby confirming the presence of the exponentially small expansion $S_1(z)$. When $x>1$, however, it is seen that $S_1(z)$ considerably exceeds the value of $B_n^{(n)}(nz)-S_0(z)$, thereby indicating its absence.
 
\begin{table}[h]
\caption{\footnotesize{Values of $B_n^{(n)}(nz)-S_0(z)$ compared with the exponentially small expansion $S_1(z)$ when $n=10$ and $z=x+i/4$}}
\begin{center}
\begin{tabular}{|c|r|r|r|}
\hline
&&&\\[-0.3cm]
\mcol{1}{|c|}{} & \mcol{1}{c|}{$x=0.60$} & \mcol{1}{c|}{$x=0.80$} & \mcol{1}{c|}{$x=0.90$}\\
\hline
&&&\\[-0.3cm]
$B_n^{(n)}(nz)-S_0(z)$ & $+0.012028+0.023460i$ & $-0.084193-0.037509i$ & $-0.089839+0.302192i$ \\
$S_1(z)$               & $+0.012023+0.023457i$ & $-0.085971-0.037707i$ & $-0.099150+0.254323i$ \\
\hline
&&&\\[-0.3cm]
\mcol{1}{|c|}{} & \mcol{1}{c|}{$x=1.10$} & \mcol{1}{c|}{$x=1.20$} & \mcol{1}{c|}{$x=1.40$}\\
\hline
&&&\\[-0.3cm]
$B_n^{(n)}(nz)-S_0(z)$ & $-0.206433-0.333096i$ & $+0.053277+0.082496i$ & $-0.018778+0.014669i$ \\
$S_1(z)$               & $+3.281489-1.068820i$ & $+6.231262-9.956311i$ & $-32.37578-94.11127i$ \\
\hline
\end{tabular}
\end{center}
\end{table}

\end{document}